\documentclass[12pt,bezier]{article}
\usepackage{amsmath}
\usepackage{amsfonts,amsthm,amssymb}
\usepackage{amsfonts}
\usepackage{graphics}
\usepackage{amssymb}
\newtheorem{lem}{Lemma}[section]
\newtheorem{thm}[lem]{Theorem}
\newtheorem{cor}[lem]{Corollary}

\theoremstyle{definition}

\begin{document}
\title{  On restricted edge-connectivity of half-transitive
multigraphs
 \footnote{The research is supported by NSFC (11401510, 11531011) and NSFXJ
(2015KL019).}}
\author{  Yingzhi Tian$^{a}$ \footnote{Corresponding author. E-mail: tianyzhxj@163.com (Y.Tian).},  Jixiang Meng$^{a}$,
Xing Chen$^{b}$
 \\
{\small $^{a}$College of Mathematics and System Sciences, Xinjiang
University,}\\
{\small Urumqi, Xinjiang, 830046, P.R.China}
\\
{\small $^{b}$Xinjiang Institute of Engineering , Urumqi, Xinjiang, 830091, P.R.China}}
\date{}

\maketitle

{\flushleft\large\bf Abstract} \ \  Let $G=(V,E)$ be a multigraph (it has multiple edges, but no loops). The edge connectivity, denoted by $\lambda(G)$, is the cardinality of a minimum edge-cut of $G$.  We call $G$ maximally edge-connected if $\lambda(G)=\delta(G)$, and $G$
super edge-connected if every minimum edge-cut is a set of edges incident with some vertex.
The restricted edge-connectivity $\lambda'(G)$ of $G$ is the minimum number of edges whose
removal disconnects $G$ into non-trivial components. If $\lambda'(G)$ achieves
the upper bound of restricted edge-connectivity, then $G$ is said to be $\lambda'$-optimal.
A bipartite multigraph is said to be half-transitive if its automorphism group is transitive on
the sets of its bipartition. In this paper, we will characterize maximally edge-connected
half-transitive multigraphs, super edge-connected half-transitive multigraphs,
and $\lambda'$-optimal half-transitive multigraphs.



\section{  Introduction  }

A graph $G$ consists of  vertex set $V$ and edge set $E$, where $E$
is a multiset of unordered pairs of (not necessarily distinct)
vertices. A $loop$ is an edge whose endpoints are the same vertex.
An edge is $multiple$ if there is another edge with the same
endvertices; otherwise it is simple. The $multiplicity$ of an edge
$e$, denoted by $\mu(e)$, is the number of multiple edges sharing
the same endvertices; the $multiplicity$ of a graph $G$, denoted by
$\mu(G)$, is the maximum multiplicity of its edges. A graph is a
$simple$ $graph$ if it has no multiple edges or loops, a
$multigraph$ if it has multiple edges, but no loops, and a
$pseudograph$ if it contains both multiple edges and loops. The
$underlying$ $graph$ of a multigraph $G$, denoted by $U(G)$, is a
simple graph obtained from $G$ by destroying all multiple edges.
It is clear that $\mu(G)=1$ if the graph $G$ is simple and contains at least one edge.

Let $G=(V,E)$ be a multigraph. The $edge$-$connectivity$ $\lambda(G)$  is the minimum size of an edge set which disconnects $G$.  Since $\lambda(G)\leq\delta(G)$, where
$\delta(G)$ is the minimum degree of $G$, a multigraph $G$ with
$\lambda(G)=\delta(G)$  is naturally said to be $maximally\
edge$-$connected$, or $\lambda$-$optimal$ for simplicity. A
multigraph $G$ is said to be $vertex$-$transitive$
if for any two vertices $u$ and $v$  in $G$, there is an automorphism $\alpha$ of $G$
such that $v=\alpha(u)$, that is, $Aut(G)$ acts transitively on $V$. A bipartite multigraph $G$
with bipartition $V_1\cup V_2$ is called $half$-$transitive$  if $Aut(G)$ acts
transitively both on $V_1$ and $V_2$. Mader \cite{Mader}
proved the following well-known result.

\begin{thm}\cite{Mader}
Every connected vertex-transitive simple graph $G$ is $\lambda$-optimal.
\end{thm}

If  $G$ is a vertex-transitive multigraph, then
$G$ is not always maximally edge-connected. A simple example is the
multigraph obtained from a 4-cycle $C_4$ by replacing each edge belonging
to a pair of opposite edges in $C_4$ with $m\ (m\geq2)$ multiple
edges.

For half-transitive simple graphs, Liang and Meng \cite{Liang}
proved the following result:

\begin{thm}\cite{Liang}
Every connected half-transitive simple graph $G$ is $\lambda$-optimal.
\end{thm}

The problem of exploring edge-connected properties stronger than the
maximally edge-connectivity for simple graphs has been widely studied. The first candidate may be the so-called $super\
edge$-$connectivity$. We can generalize this definition to multigraphs.
A multigraph $G$ is said to be $super\ edge$-$connected$, in short, $super$-$\lambda$,
if each of its
minimum edge-cuts isolates a vertex, that is, every minimum
edge-cut is a set of edges incident with a certain vertex in $G$. By
the definitions, a super-$\lambda$ multigraph must be a
$\lambda$-optimal multigraph. However, the converse is not true. For
example, the multigraph obtained from $K_m\times K_2$ by replacing every edge with a pair of multiple edges is $\lambda$-optimal but not
super-$\lambda$ since the set of edges between the two copies of the multi-subgraph obtained from $K_m$ by replacing every edge with a pair of multiple edges is a minimum edge-cut which does not isolate any vertex.

The concept of super-$\lambda$  was originally introduced by Bauer
et al. \cite{Bauer}, where combinatorial optimization problems
in the design of reliable probabilistic simple graphs were investigated. The
following theorem is a nice result of Tindell
\cite{Tindell1}, who characterized super edge-connected
vertex-transitive simple graphs.

\begin{thm} \cite{Tindell1}
A connected vertex-transitive simple graph $G$ which is neither a
cycle nor a complete graph  is super-$\lambda$ if and only if it
contains no clique $K_k$, where k is the degree of $G$.
\end{thm}

For further study, Esfahanian and Hakimi \cite{Esfahanian2}
introduced the concept of restricted edge-connectivity for simple
graphs. The concept of restricted edge-connectivity is one kind of
conditional edge-connectivity proposed by Harary in \cite{Harary},
and has been successfully applied in the further study of tolerance
and reliability of networks, see [2,6,8,11-12,18,20-23]. Let $F$ be
a set of edges in $G$. Call $F$ a $restricted\ edge$-$cut$ if $G-F$
is disconnected and contains no isolated vertices. The minimum
cardinality over all restricted edge-cuts is called $restricted\
edge$-$connectivity$ of $G$, and denoted by $\lambda'(G)$. It was
shown by Wang and Li \cite{Wangm} that the larger $\lambda'(G)$ is,
the more reliable the network is. In \cite{Esfahanian2}, it was
proved that if a connected simple graph $G$ of order $|V(G)|\geq4$
is not a star $K_{1,n-1}$, then $\lambda'(G)$ is well-defined and
$\lambda'(G)\leq \xi(G)$, where $\xi(G)=$min$\{d(u)+d(v)-2:uv\in
E(G)\}$ is the minimum edge degree of $G$. A simple graph $G$ with
$\lambda'(G)=\xi(G)$ is called a $\lambda'$-$optimal\ graph$. It
should be pointed out that if $\delta(G)\geq3$, then a
$\lambda'$-optimal simple graph must be super-$\lambda$. In fact, a
graph $G$ is super-$\lambda$ if and only if
$\lambda(G)<\lambda'(G)$, see \cite{Li1}. Thus, the concepts of
$\lambda$-optimal graphs, super-$\lambda$ graphs and
$\lambda'$-optimal graphs describe reliable interconnection
structures for graphs at different levels.

In \cite{Meng1}, Meng studied  the parameter $\lambda'$
for connected vertex-transitive simple graphs. The main result  may be restated as follows:

\begin{thm} \cite{Meng1}
Let $G$ be a $k$-regular connected vertex-transitive simple graph
which is neither a cycle nor a complete graph. Then $G$ is not
$\lambda'$-optimal if and only if it contains a $(k-1)$-regular
subgraph $H$ satisfying $k\leq |V(H)|\leq 2k-3$.
\end{thm}

The authors in \cite{Tian2} proved the following result.

\begin{thm} \cite{Tian2}
Let $G=(V_1\cup V_2, E)$ be a connected half-transitive simple graph with
$n=|V(G)|\geq4$ and $G\ncong K_{1,n-1}$. Then $G$ is
$\lambda'$-optimal.
\end{thm}

Since a graph $G$ is super-$\lambda$ if and only if
$\lambda(G)<\lambda'(G)$, Theorem 1.5 implies the following corollary.

\begin{cor}
The only connected half-transitive simple graphs which are not
super-$\lambda$ are cycles $C_n(n\geq4)$.
\end{cor}

We can naturally generalize the concept of restricted edge-connectivity to
multigraphs. The $restricted$ $edge$-$connectivity$ $\lambda'(G)$ of
a multigraph $G$ is the minimum number of edges whose removal
disconnects $G$ into non-trivial components. Similarly, define the
minimum edge degree of $G$ as $\xi(G)=min\{\xi(e)=d(u)+d(v)-2\mu(e):
e=uv\in E(G)\}$, where $\xi(e)=d(u)+d(v)-2\mu(e)$ is the edge degree
of the edge $e=uv$ in $G$. But the inequality
$\lambda'(G)\leq \xi(G)$ is not always correct. For example, the
restricted edge-connectivity of the multigraph $G$ in Fig.1 is 6, but
$\xi(G)=4$.

In \cite{Tian3}, we gave sufficient and necessary conditions for vertex-transitive multigraphs
to be maximally edge-connected, super edge-connected and $\lambda'$-optimal.
In the following, we will study maximally edge-connected half-transitive multigraphs,
super edge-connected half-transitive multigraphs, and $\lambda'$-optimal half-transitive multigraphs.

\begin{center}
\scalebox{0.15}{\includegraphics{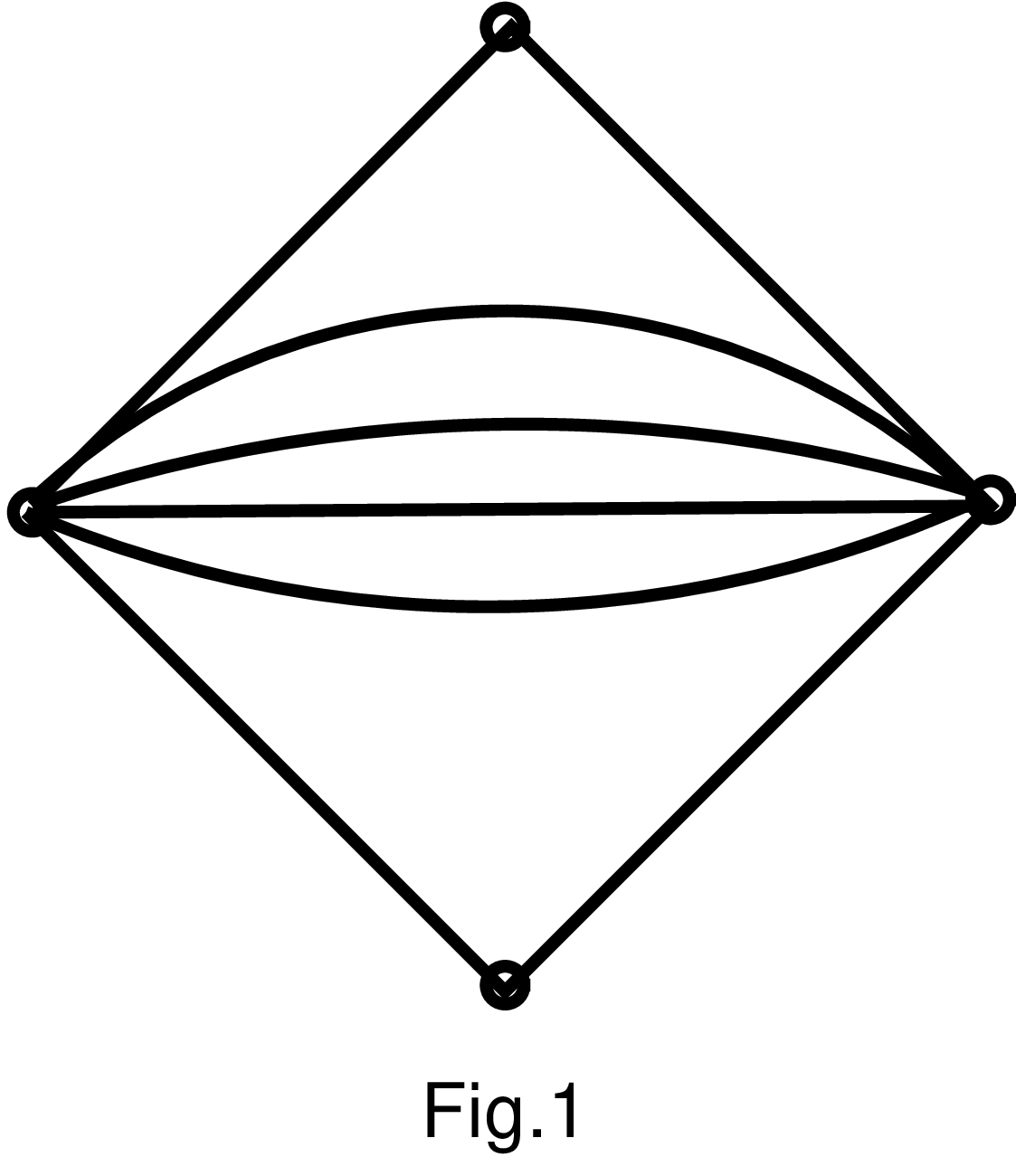}}
\end{center}

\section{Preliminaries}

Let $G=(V,E)$ be a multigraph. For two disjoint non-empty subsets $A$ and
$B$ of $V$, let $[A,B]=\{e=uv\in E:u\in A \ \mbox{and} \ v\in B\}$.
For the sake of convenience, we write $u$ for the single vertex set
$\{u\}$. If $\overline{A}=V\backslash A$, then we write $N(A)$
for $[A,\overline{A}]$ and $d(A)$ for $|N(A)|$. Thus $d(u)$ is just the degree of $u$ in $G$.
Denote by $G[A]$ the subgraph of $G$ induced by $A$.

An edge-cut $F$ of $G$ is called a $\lambda$-cut if
$|F|=\lambda(G)$. It is easy to see that for any $\lambda$-cut $F$,
$G-F$ has exactly two components. If $N(A)$ is a
$\lambda$-cut of $G$, then $A$ is called a $\lambda$-$fragment$ of
$G$. It is clear that if $A$ is a $\lambda$-fragment of $G$, then so
is $\overline{A}$. Let $r(G)$=min$\{|A|$: $A$ is a
$\lambda$-fragment of $G\}$. Obviously, $1\leq r(G)\leq
\frac{1}{2}|V|$. A $\lambda$-fragment $B$ is called a
$\lambda$-$atom$ of $G$ if $|B|=r(G)$. A $\lambda$-fragment $C$ is
called a $strict$ $\lambda$-$fragment$ if $2\leq |C| \leq |V(G)|-2$.
If $G$ contains strict $\lambda$-fragments, then the ones with
smallest cardinality are called $\lambda$-$superatoms$.

Similarly, we can give the definition of $\lambda'$-$atom$.
A restricted edge-cut $F$ of $G$ is called a $\lambda'$-cut if
$|F|=\lambda'(G)$. For any $\lambda'$-cut
$F$, $G-F$ has exactly two components. Let $A$ be a proper
subset of $V$. If $N(A)$ is a $\lambda'$-cut of $G$, then $A$
is called a $\lambda'$-$fragment$ of $G$. It is clear that if $A$ is
a  $\lambda'$-fragment of $G$, then so is $\overline{A}$. Let
$r'(G)$=min$\{|A|$: $A$ is a $\lambda'$-fragment of $G\}$.
Obviously, $2\leq r'(G)\leq \frac{1}{2}|V|$. A $\lambda'$-fragment
$B$ is called a $\lambda'$-$atom$ of $G$ if $|B|=r'(G)$.

For a multigraph $G$, the inequality $\lambda'(G)\leq \xi(G)$ is not always correct.
But if $G$ is a $k$-regular multigraph, we proved the following result.

\begin{lem}\label{2-lem-1} \cite{Tian3}
Let $G$ be a connected $k$-regular multigraph. Then $\lambda'(G)$ is well-defined and
$\lambda'(G)\leq\xi(G)$ if  $|V(G)|\geq4$.
\end{lem}

We call a bipartite multigraph $G$ with bipartition $V_1\cup V_2$ $semi$-$regular$ if
each vertex in $V_1$ has the same degree $d_1$ and each vertex in $V_2$ has the same degree $d_2$ in $G$.
For semi-regular bipartite multigraphs,  a similar result can be obtained.

\begin{lem}\label{2-lem-1}
Let $G$ be a connected semi-regular bipartite multigraph with bipartition $V_1\cup V_2$.
Then $\lambda'(G)$ is well-defined and $\lambda'(G)\leq\xi(G)$ if  $|V(G)|\geq4$ and $U(G)\ncong K_{1,n-1}$.
\end{lem}

\noindent{\bf Proof.} Assume each vertex in $V_1$ has degree $d_1$ and each vertex in $V_2$ has
 degree $d_2$ in $G$. Assume, without loss of generality, that $d_1\leq d_2$. Let $e=uv$ be an edge
 such that $\xi(e)=\xi(G)$, where $u\in V_1$ and $v\in V_2$. If $G-\{u,v\}$
contains a non-trivial component, say $C$, then $N(V(C))$ is a restricted edge-cut
and $|N(V(C))|\leq |N(\{u,v\})|=\xi(e)=\xi(G)$. Thus assume that $G-\{u,v\}$ only contains isolated vertices.
 If there is a vertex $w$ other than $v$ in $V_2$, then since $|V\setminus\{u,v\}|=|V(G)|-2\geq2$, we obtain a contradiction
$d_1+d_2\leq|N(V\setminus\{u,v\})|=|N(\{u,v\})|=\xi(e)=d_1+d_2-2\mu(e)<d_1+d_2$. Therefore, $V_2=\{v\}$ and  $U(G)\cong K_{1,n-1}$, also a contradiction. $\square$

Because of Lemma 2.1 and Lemma 2.2, we call a regular multigraph (or a semi-regular bipartite multigraph)
$G$ $\lambda'$-$optimal$ if $\lambda'(G)=\xi(G)$. Since each vertex-transitive
multigraph is regular and each half-transitive multigraph is semi-regular,
thus a vertex-transitive multigraph (or  a half-transitive multigraph)
$G$ is $\lambda'$-$optimal$ if $\lambda'(G)=\xi(G)$.

Recall that an $imprimitive$ $block$ for a permutation group $\Phi$
on a set $T$ is a proper, non-trivial subset $A$ of $T$ such that
for every $\varphi\in\Phi$  either  $\varphi(A)=A$ or
$\varphi(A)\cap A=\O$. A subset $A$ of $V(G)$ is called an
$imprimitive$ $block$ for $G$ if it is an imprimitive block for the
automorphism group $Aut(G)$ on $V(G)$. The following theorem
shows the importance of imprimitive blocks:

\begin{thm}\label{pro-1} \cite{Tindell2} Let $G=(V,E)$ be a connected
simple graph and $A$ be an imprimitive block for $G$.  If $G$ is
vertex-transitive, then $G[A]$ is also vertex-transitive.
\end{thm}

By a similar argument as Theorem 2.3, we can obtain the following result for half-transitive multigraphs.

\begin{lem} Let $G$ be a connected bipartite multigraph with bipartition $V_1\cup V_2$.
Assume $A$ is an imprimitive block for $G$ such that $A\cap V_1\neq\O$ and $A\cap V_2\neq\O$.
 If $G$ is half-transitive, then $G[A]$ is also half-transitive.
\end{lem}

\noindent{\bf Proof.} Since $G$ is half-transitive, for any two vertices $u,v\in A\cap V_i$ ($i\in \{1,2\}$),
 there is  $\alpha\in Aut(G)$ such that $\alpha(u)=v$. Because $\alpha(A)\cap A\neq\O$,
 we have $\alpha(A)=A$ by $A$ is an imprimitive block for $G$.
 Thus the restriction of $\alpha$ to $A$ is an automorphism of $G[A]$,
which maps $u$ to $v$. It follows  that $G[A]$ is a  half-transitive multigraph. $\square$

\begin{lem}
Let $G$ be a connected half-transitive multigraph with bipartition $V_1\cup V_2$ and $A$ be an imprimitive block for $G$ with $A_1=A\cap V_1\neq\emptyset$ and $A_2=A\cap V_2\neq\emptyset$. Assume each vertex in $V_1$ has degree $d_1$ and each vertex in $V_2$ has degree $d_2$ in $G$, and each vertex in $A_1$ has degree $d_1'$ and each vertex in $A_2$ has degree $d_2'$ in $G[A]$. Then $d_1'<d_1$ and $d_2'<d_2$.
\end{lem}

\noindent{\bf Proof.} Since $G$ is half-transitive, for $u_i\in A_i$ and $v_i\in V_i\backslash A_i$ ($i\in\{1,2\}$), there exists an automorphism $\alpha\in Aut(G)$ such that $\alpha(u_i)=v_i$. Because $A$ is an imprimitive block for $G$ and $\alpha(u_i)\notin A$, we have $\alpha(A)\cap A=\emptyset$. Thus there exist $\alpha_1,\alpha_2,\cdots,\alpha_p\in Aut(G)$
satisfying $V(G)=\cup_{i=1}^{p}\alpha_i(A)$ and $\alpha_i(A)\cap\alpha_j(A)=\emptyset$ for $1\leq i\neq j\leq p$. Since $G$ is connected and $G[\alpha_i(A)]\cong G[A]$ for $1\leq i\leq p$, we can verify that $d_1'<d_1$ and $d_2'<d_2$. $\square$

\section{Maximally edge-connected half-transitive multigraphs}

In \cite{Mader}, Mader proved that  any two distinct $\lambda$-atoms of a simple graph are disjoint.
 For multigraphs, this property still holds.

\begin{lem}\label{2-lem-1}
Let $G$ be a connected multigraph. Then any two distinct $\lambda$-atoms of $G$ are disjoint.
\end{lem}

\noindent{\bf Proof.} Suppose to the contrary that there are two distinct $\lambda$-atoms $A$ and $B$
with $A\cap B\neq \O$. We have $V(G)\backslash (A\cup B)\neq \O$ by $|A|\leq |V(G)|/2$ and $|B|\leq |V(G)|/2$.
Then $N(A\cap B)$ and  $N(A\cup B)$ are edge-cuts of $G$, thus $d(A\cap B)=|N(A\cap B)|\geq\lambda(G)$ and
  $d(A\cup B)=|N(A\cup B)|\geq\lambda(G)$. From the following well-known submodular inequality (see \cite{Tindell2}),
$$2\lambda(G)\leq d(A\cup B)+d(A\cap B)\leq d(A)+d(B)=2\lambda(G),$$
we conclude that  both  $d(A\cap B)=\lambda(G)$ and  $d(A\cup B)=\lambda(G)$ hold. Keep in mind, $d(A\cap B)=\lambda(G)$ implies that $N(A\cap B)$ is a minimum edge-cut and thus both $A\cap B$ and $V\backslash (A\cap B)$ are connected. Therefore, $A\cap B$ is a $\lambda$-fragment with $|A\cap B|<|A|$, which contradicts to $A$
is a $\lambda$-atom of $G$.
$\square$

\begin{thm}\label{2-lem-1}
Let $G$ be a connected half-transitive multigraph with bipartition $V_1\cup V_2$.
 Assume each vertex in $V_1$ has degree $d_1$ and each vertex in $V_2$ has degree $d_2$ in $G$.
 Then $G$ is not maximally edge-connected if and only if there is a proper induced  connected half-transitive
multi-subgraph $H$ of $G$ such that
$$d_1'<d_1,d_2'<d_2\ and\ |A_1|(d_1-d_1')+|A_2|(d_2-d_2')\leq min\{d_1,d_2\}-1,$$ where
$A_1=V_1\cap V(H)$, $A_2=V_2\cap V(H)$, $d_1'$ is the degree of each
vertex of $A_1$ and $d_2'$ is the degree of each vertex of
$A_2$ in $H$.
\end{thm}

\noindent{\bf Proof.} Assume, without loss of generality, that
$d_1\leq d_2$. If $G$ is not maximally edge-connected, then
$\lambda(G)\leq d_1-1$. Let $A$ be a $\lambda$-atom of $G$ and
$H=G[A]$. By Lemma 3.1, we know that $A$ is an imprimitive block for $G$.
Thus $H$ is a connected half-transitive multigraph by Lemma 2.4.
Assume each vertex in $A\cap V_1$ has degree $d_1'$ and each
vertex in $A\cap V_2$ has degree $d_2'$ in $H$. Then $|A\cap
V_1|(d_1-d_1')+|A\cap V_2|(d_2-d_2')=d(A)=\lambda(G)\leq d_1-1$.
By Lemma 2.5, $d_1'<d_1$ and $d_2'<d_2$.

Now we prove the sufficiency. Assume $G$ contains a proper induced
connected half-transitive multi-subgraph $H$ such that
$d_1'<d_1,d_2'<d_2\ and\ |A_1|(d_1-d_1')+|A_2|(d_2-d_2')\leq min\{d_1,d_2\}-1$, then
$\lambda(G)\leq d(V(H))=|A_1|(d_1-d_1')+|A_2|(d_2-d_2')\leq
min\{d_1,d_2\}-1$, that is, $G$ is not maximally edge-connected.
$\square$

%

\section{Super edge-connected half-transitive multigraphs}

In \cite{Tindell2}, Tindell studied the intersection property of $\lambda$-superatoms of
vertex-transitive simple graphs. For half-transitive multigraphs,  we have the following lemma.

\begin{lem}\label{2-lem-1}
Let $G$ be a connected half-transitive multigraph
with bipartition $V_1\cup V_2$. Assume $G$ is not super
edge-connected, $A$ and $B$ are two distinct $\lambda$-superatoms.
If $|A|=|B|\geq3$,  then $A\cap B=\O$.
\end{lem}

\noindent{\bf Proof.} Assume each vertex in $V_1$ has degree $d_1$ and each vertex in $V_2$ has degree $d_2$ in $G$.
Without loss of generality, assume that $d_1\leq d_2$.
If $A\cap B\neq\O$, then by a similar
argument as the proof of Lemma 3.1, we can conclude that $d(A\cap
B)=d(A\cup B)=\lambda(G)$. Since $|A|=|B|$ and $A\neq B$, we know that $|A\cap B|\leq |V(G)|-2$. Hence, if $|A\cap B|\geq 2$, then 
it is a strict $\lambda$-fragment strictly contained in $A$ which contradicts to $A$ being a $\lambda$-superatom (Because $d(A\cap B)=\lambda(G)$ implies that $N(A\cap B)$ is a minimum edge-cut and thus both $A\cap B$ and $V\backslash (A\cap B)$ are connected). Therefore $|A\cap B|=1$.

Let $C=V(G)\setminus B$. Then $|A\cap C|=|A\setminus (A\cap
B)|\geq2$, and $A$, $V(G)\setminus A$, $C$  and $V(G)\setminus C$
are all strict $\lambda$-fragments. By a similar argument as above we can
deduce that $A \cap C$ is a strict $\lambda$-fragment with $|A\cap C|<|A|$, which is
impossible.
$\square$

\begin{thm}\label{2-lem-1}
Let $G$ be a connected half-transitive multigraph
with bipartition $V_1\cup V_2$. Assume each vertex in $V_1$ has
degree $d_1$, each vertex in $V_2$ has degree $d_2$ in $G$  and
$|V(G)|\geq2\ min\{d_1,d_2\}+2$. Then $G$ is not super
edge-connected if and only if there is a proper induced connected
half-transitive multi-subgraph $H$ of $G$ such that
$$d_1'<d_1,d_2'<d_2\ and\ |A_1|(d_1-d_1')+|A_2|(d_2-d_2')\leq min\{d_1,d_2\},$$ where $A_1=V_1\cap
V(H)$, $A_2=V_2\cap V(H)$, $d_1'$ is the degree of each vertex of
$A_1$ and $d_2'$ is the degree of each vertex of $A_2$ in
$H$.
\end{thm}

\noindent{\bf Proof.} Assume, without loss of generality, that
$d_1\leq d_2$. If $G$ is not super edge-connected, then $G$ contains
$\lambda$-superatoms. Let $A$ be a $\lambda$-superatom of $G$ and
$H=G[A]$. If $|A|=2$, then $H$ is isomorphic to a multigraph which
contains two vertices and $t$ edges between the two vertices. Thus
$H$ is an induced $t$-regular connected half-transitive
multi-subgraph of $G$. Therefore $|A\cap V_1|(d_1-t)+|A\cap
V_2|(d_2-t)=d(A)=\lambda(G)\leq d_1$. Since $G$ is both connected and  half-transitive, we can verify that $t<d_1$. In the following, we assume that
$|A|\geq3$.

Lemma 4.1 impies that $A$ is an imprimitive block for $G$. Thus $H$
is a connected half-transitive multigraph by Lemma 2.4. Assume each
vertex in $A\cap V_1$ has degree $d_1'$ and each vertex in
$A\cap V_2$ has degree $d_2'$ in $H$. Thus $|A\cap
V_1|(d_1-d_1')+|A\cap V_2|(d_2-d_2')=d(A)=\lambda(G)\leq d_1$. By Lemma 2.5, $d_1'<d_1$ and $d_2'<d_2$.

Now we prove the sufficiency. If $\lambda<min\{d_1,d_2\}$, then $G$ is not super edge-connected. Therefore, we only need to consider the case when $\lambda=min\{d_1,d_2\}$. Assume $G$ contains a proper induced
connected half-transitive multi-subgraph $H$ such that
$d_1'<d_1,d_2'<d_2\ and\ |A_1|(d_1-d_1')+|A_2|(d_2-d_2')\leq min\{d_1,d_2\}$, then
$d(V(H))=|A_1|(d_1-d_1')+|A_2|(d_2-d_2')\leq min\{d_1,d_2\}$. If
$G-V(H)$ contains no isolated vertices, then $V(H)$ is a strict
$\lambda$-fragment. Thus $G$ is not super edge-connected. Assume
$G-V(H)$ contains an isolated vertex $w$. Then $N(w)=N(V(H))$.
Since $|A_1|\leq\ min\{d_1,d_2\}$ and $|A_2|\leq\ min\{d_1,d_2\}$ by
$d_1'<d_1,d_2'<d_2\ and\ |A_1|(d_1-d_1')+|A_2|(d_2-d_2')\leq min\{d_1,d_2\}$, we see that $G$
is not connected by $|V(G)|\geq2\ min\{d_1,d_2\}+2$, a
contradiction. $\square$

\section{$\lambda'$-optimal half-transitive multigraphs}

In \cite{Xu}, the authors proved the following fundamental result for
studying the restricted edge-connectivity of simple graphs.

\begin{thm} \cite{Xu}
Let $G=(V,E)$ be a connected simple graph with at least four
vertices and $G\ncong K_{1,n-1}$. If $G$ is not $\lambda'$-optimal,
then any two distinct $\lambda'$-atoms of $G$ are disjoint.
\end{thm}

For multigraphs, we cannot obtain a similar result as in Theorem
5.1. But for half-transitive multigraphs, the similar result holds.

\begin{lem}\label{2-lem-1}
Let $G$ be a connected multigraph with $\delta(G)\geq2\mu(G)$. If
$G$ contains a $\lambda'$-atom $A$ with $|A|\geq3$, then each vertex
in $A$ has at least two neighbors in $A$.
\end{lem}

\noindent{\bf Proof.} By contradiction, assume there is a vertex
$u\in A$ such that $u$ contains only one neighbor in $A$. Let $v$ be
the only neighbor of $u$ in $A$. Set $A'=A\backslash\{u\}$. Then
both $G[A']$ and $G[\overline{A'}]$ are connected. We have
$|A'|\geq2$ by $|A|\geq3$. Clearly,
$|\overline{A'}|=|\overline{A}|+1\geq4$. Thus $[A',\overline{A'}]$
is a restricted edge-cut.  Since $\delta(G)\geq2\mu(G)$, we have
$$\lambda'(G)\leq
|[A',\overline{A'}]|=|[A,\overline{A}]|+\mu(uv)-(d(u)-\mu(uv))\leq
|[A,\overline{A}]|=\lambda'(G).$$ It follows that $A'$ is a
$\lambda'$-fragment with $|A'|<|A|$, which contradicts to $A$ is a
$\lambda'$-atom. $\square$

The proof of Lemma 5.3 is inspired by [13, Lemma 4.2].

\begin{lem}\label{2-lem-1}
Let $G$ be a connected half-transitive multigraph
with bipartition $V_1\cup V_2$ and $\delta(G)\geq2\mu(G)$. Assume
$G$ is not $\lambda'$-optimal, $A$ and $B$ are two distinct
$\lambda'$-atoms. Then $|A|=|B|\geq3$ and $A\cap B=\O$.
\end{lem}

\noindent{\bf Proof.} Assume each vertex in $V_1$ has degree $d_1$
and each vertex in $V_2$ has degree $d_2$ in $G$. Without loss of
generality, assume that $d_1\leq d_2$.

If $|A|=2$, then $\lambda'(G)=d(A)=d_1+d_2-2\mu(uv)\geq\xi(G)$
(where $A=\{u,v\}$), which contradicts that $G$ is not
$\lambda'$-optimal. Thus $|A|\geq3$.

Suppose to the contrary that $A\cap B\neq \O$. Set $C=A\cap B$,
$A_1=A\cap\overline{B}$, $B_1=B\cap \overline{A}$ and
$D=\overline{A}\cap \overline{B}=\overline{A \cup B}$. In the
following, we will derive a contradiction by a series of claims.

Clearly, one of the following two inequalities must hold:
$$|[A_1,C]|\leq|[B_1,C]|+|[C,D]|,\eqno(1)$$
$$|[B_1,C]|\leq|[A_1,C]|+|[C,D]|.\eqno(2)$$
In the following, we always assume, without loss of generality, that
inequality (1) holds.

\noindent{\bf Claim 1.} $A_1$ satisfies one of the following two
conditions: ($i$) $A_1=\{v_{21}\} (v_{21}\in V_2)$ and
$d_1>2\mu(G)$, or ($ii$) $A_1=\{v_{11},\cdots,v_{1m}\}(v_{1i}\in
V_1$ for $1\leq i\leq m$) and $d_2>(m-1)d_1+2\mu(G)$.

It follows from inequality (1) that
$$d(A_1)=|[A_1,D]|+|[A_1,C]|+|[A_1,B_1]|\leq d(A)=\lambda'(G).$$

Assume $G[A_1]$ has a component $\widetilde{G}$ with
$|V(\widetilde{G})|\geq2$. Set $F=V(\widetilde{G})$. Since $G[B]$
and $G[\overline{A}]$ are both connected, and $B\cap
\overline{A}\neq\O$, we see that $G[\overline{A_1}]$ is connected.
Furthermore, since $G$ is connected, every component of $G[A_1]$ is
joined to $G[\overline{A_1}]$, and thus $G[\overline{F}]$ is
connected. So $[F,\overline{F}]$ is a restricted edge-cut with
$d(F)\leq\lambda'(G)$. Because $A$ is a $\lambda'$-atom and $F$ is
a proper subset of $A$, we obtain $d(F)>d(A)=\lambda'(G)$, a
contradiction. Thus, each component in $G[A_1]$ is an isolated
vertex. By $d(A_1)\leq \lambda'(G)< d_1+d_2-2\mu(G)$, we can derive
that $A_1$ satisfies one of the following two conditions: ($i$)
$A_1=\{v_{21}\} (v_{21}\in V_2)$ and $d_1>2\mu(G)$, or ($ii$)
$A_1=\{v_{11},\cdots,v_{1m}\}(v_{1i}\in V_1$ for $1\leq i\leq m$)
and $d_2>(m-1)d_1+2\mu(G)$.

\noindent{\bf Claim 2.} $C\nsubseteq V_1$ and $C\nsubseteq V_2$.

By contradiction. Suppose $C\subseteq V_1$. Then $G[C]$ is an
independent set. Since we have assumed that
$|[A_1,C]|\leq|[C,B_1]|+|[C,D]|$, there exists a vertex $v$ in $C$
such that
$$|[v,A_1]|\leq|[v,D]| + |[v,B_1]|.\eqno(3)$$
Set $F=A\setminus \{v\}$, then
$$d(F)=d(A)-|[v,D]|-|[v,B_1]|+|[v,A_1]|\leq d(A)=\lambda'(G).$$
Since $G[A]$ is connected and $C$ is an independent set, we have
$|[v,A_1]|\geq1$. It follows from inequality (3) that
$|[v,\overline{A}]|\geq1$. So, $G[\overline{F}]$ is connected. We
claim that each component in $G[F]$ has at least 2 vertices. Indeed, if there is an isolated vertex $u$ in $G[F]$, then $v$ is the
only vertex adjacent to $u$ in $G[A]$, which contradicts to Lemma
5.2.  Now, similarly as in the proof of Claim 1, a contradiction
arises, since $F$ contains a smaller $\lambda'$-fragment than
$A$.  $C\nsubseteq V_2$ can be proved similarly.

\noindent{\bf Claim 3.} $d(D)<\lambda'(G)$ and $D$ is an independent
set contained in $V_1$.

By Claim 2, $|C|\geq2$. We claim that $d(C)>\lambda'(G)$. In fact,
if $G[C]$ contains a component of order at least 2, then similar to
the proof of Claim 1, we can show that $[C,\overline{C}]$ contains a
restricted edge-cut, and thus $d(C)>\lambda'(G)$. Otherwise,
we assume that each component in $G[C]$ is  an
isolated vertex. Since not all
vertices in $C$ are from the same bipartition, there must be at
least one vertex in $V_2$. From $|C|\geq2$, we have $d(C)\geq
d_2+d_1>\xi(G)\geq\lambda'(X)$. Thus, we have that $d(C)>\lambda'(G)$.

From the well-known submodular  inequality (see \cite{Tindell2}), we
have
$$d(C)+d(D)\leq d(A)+d(B)=2\lambda'(G).\eqno(4)$$
By (4) and $d(C)>\lambda'(G)$, we obtain $d(D)<\lambda'(G)$.
Applying a similar argument as above, we can show that $D$ is an
independent set contained in $V_1$.

Since $|A_1|+|C|=|A|\leq|\overline{A}|=|B_1|+|D|$ and $|A_1|=|B_1|$, we have
$|D|\geq|C|$. Let $s=|D|$. Then $s\geq|C|\geq2$ and
$$d(D)=sd_1.\eqno(5)$$

Denote by $e_1$ the number of edges in $G[\overline{C}]$. Clearly,
$$d(C)=d(\overline{C})=\sum_{v\in \overline{C}}d(v)-2e_1.\eqno(6)$$
Since $G[\overline{B}]$ is connected and $D$ is an independent set
contained in $V_1$, Claim 1 ($ii$) can not hold. Thus, Claim 1 ($i$)
is true. This implies $|A_1|=|B_1|=1$. Because $G[\overline{A}]$ is connected and $D$ is an independent set contained in $V_1$, we know $B_1\subseteq V_2$. Since $G$ is a bipartite multigraph, we have
$$e_1\leq2s\mu(G).\eqno(7)$$
Combining this with (4), (5) and (6), we see that
$$2d_1+2d_2-4\mu(G)-sd_1>2\lambda'(G)-d(D)\geq d(C)\geq sd_1+2d_2-4s\mu(G).$$
This implies $d_1<2\mu(G)$, contradicting to the assumption that
$d_1\geq2\mu(G)$.  $\square$

\begin{thm}\label{2-lem-1}
Let $G$ be a connected half-transitive multigraph
with bipartition $V_1\cup V_2$ and $\delta(G)\geq2\mu(G)$. Assume
each vertex in $V_1$ has degree $d_1$, each vertex in $V_2$ has
degree $d_2$ in $G$, $|V_1|\geq\xi(G)$ and $|V_2|\geq\xi(G)$. Then $G$ is
not $\lambda'$-optimal if and only if there is a proper induced
connected half-transitive multi-subgraph $H$ of $G$ such that
$$d_1'<d_1,d_2'<d_2\ and\ |A_1|(d_1-d_1')+|A_2|(d_2-d_2')\leq\xi(G)-1,$$ where $A_1=V_1\cap
V(H)$, $A_2=V_2\cap V(H)$, $d_1'$ is the degree of each vertex of
$A_1$ and $d_2'$ is the degree of each vertex of $A_2$ in
$H$.
\end{thm}

\noindent{\bf Proof.} Assume $G$ is not $\lambda'$-optimal. By Lemma 2.2, $G$ contains
$\lambda'$-atoms. Let $A$ be a $\lambda'$-atom of $G$  and $H=G[A]$.
By Lemma 5.3, we have $|A|\geq3$ and $A$ is an imprimitive block for
$G$. Thus $H$ is a connected half-transitive multigraph by Lemma
2.4. Assume each vertex in $A\cap V_1$ has degree $d_1'$ and
each vertex in $A\cap V_2$ has degree $d_2'$ in $H$. Then $|A\cap
V_1|(d_1-d_1')+|A\cap V_2|(d_2-d_2')=d(A)=\lambda'(G)\leq\xi(G)-1$.
By Lemma 2.5, $d_1'<d_1$ and $d_2'<d_2$.

Now we prove the sufficiency. Assume $G$ contains a proper induced
connected half-transitive multi-subgraph $H$ such that
$d_1'<d_1,d_2'<d_2\ and\ |A_1|(d_1-d_1')+|A_2|(d_2-d_2')\leq\xi(G)-1$, then
$d(V(H))=|A_1|(d_1-d_1')+|A_2|(d_2-d_2')\leq\xi(G)-1$,
$|A_1|\leq\xi(G)-1$ and $|A_2|\leq\xi(G)-1$. If $G-V(H)$
contains a non-trivial component, say $B$, then $[B,\overline{B}]$ is a restricted edge-cut and $d(B)\leq d(V(H))\leq \xi(G)-1$. Thus $G$ is
not $\lambda'$-optimal. Now we assume that each component of
$G-V(H)$ is an isolated vertex, then $d(V(\overline{H}))\geq
d_1+d_2>\xi(G)$ by $|V_1|\geq\xi(G)$ and $|V_2|\geq\xi(G)$. On the
other hand, $d(V(\overline{H}))=d(V(H))\leq \xi(G)-1$, 
a contradiction. $\square$

\end{document}